\documentclass[journal]{IEEEtran}

\usepackage{graphics} 
\usepackage{epsfig} 
\usepackage{mathptmx} 
\usepackage{times} 
\usepackage{amsmath} 
\usepackage{amssymb}  

\title{\LARGE \bf
The Role of Singular Control in Frictionless Atom Cooling in a Harmonic Trapping Potential
}

\author{Dionisis Stefanatos and Jr-Shin Li
\thanks{This work was supported by the Air Force Office of Scientific Research under Young Investigator Award FA9550-10-1-0146}
\thanks{D. Stefanatos and J.-S. Li are with the Department of Electrical and Systems Engineering, Washington University,
        St. Louis, MO 63130, USA
        {\tt\small dionisis@seas.wustl.edu, jsli@seas.wustl.edu}}%
}

\begin{document}

\maketitle
\thispagestyle{empty}
\pagestyle{empty}

\begin{abstract}

In this article we study the frictionless cooling of atoms trapped in a harmonic potential, while minimizing the transient energy of the system. We show that in the case of unbounded control, this goal is achieved by a singular control, which is also the time-minimal solution for a ``dual" problem, where the energy is held fixed. In addition, we examine briefly how the solution is modified when there are bounds on the control. The results presented here have a broad range of applications, from the cooling of a Bose-Einstein condensate confined in a harmonic trap to adiabatic quantum computing and finite time thermodynamic processes.

\end{abstract}

\section{INTRODUCTION}

Frictionless atom cooling in a harmonic trapping potential is defined as the problem of
changing the harmonic frequency of the trap to some lower final value, while
keeping the populations of the initial and final levels invariant, thus
without generating friction and heating. Achieving this goal in minimum time
has many important potential applications. For example, it can be used to
reach extremely low temperatures inaccessible by standard cooling techniques
\cite{Leanhardt03}, to reduce the velocity dispersion and collisional shifts
for spectroscopy and atomic clocks \cite{Bize05}, and in adiabatic quantum
computation \cite{Aharonov07}. It is also closely related to the problem of
moving in minimum time a system between two thermal states, as for example in
the transition from graphite to diamond \cite{Salamon09}.

The quest for optimal controls that minimize the necessary time for the cooling process has an interesting history.
It was initially proved that minimum transfer time can be achieved with bang-bang real
frequency controls \cite{Salamon09}. Later, it was
shown that when the restriction for real frequencies is relaxed, allowing the
trap to become an expulsive parabolic potential at some time intervals,
shorter transfer times can be obtained \cite{Chen10a}. Based on these previous works, we formulated frictionless atom
cooling as a minimum-time optimal control problem, permitting the frequency to
take real and imaginary values in specified ranges \cite{Stefanatos10}. We showed that the optimal
solution has again a bang-bang form and used this fact to obtain estimates of
the minimum transfer times for various numbers of switchings \cite{Stefanatos10}. We subsequently fully solved the corresponding
time-optimal control problem and obtained the optimal synthesis \cite{Stefanatos11}.

Recently, it was pointed out that the energy ``cost" of the cooling process, more precisely the transient excitation energy, can impose limits to the possible speed-up \cite{Chen10b}.
For example, for a trap which is harmonic near the ground state, but not for higher energies, large transient energies imply perturbing effects of anharmonicities and thus undesired excitations of the final state. The problem of minimizing the transient energy for a fixed transfer time and with unlimited controls was considered in \cite{Chen10b}, and an interesting conjugate relation between these two quantities was revealed.

In the present article we examine this problem from a control-theoretic point of view. We show that the transient energy is minimized by a singular control. We also show that this singular control is the time-minimal solution for a ``dual" problem, where the energy is held fixed, elucidating thus the conjugate relation between transient energy and minimum time. Finally, we examine briefly how the solution is modified when there are bounds on the control. The results presented here come to meet several examples of singular solutions for optimal control problems on spin systems existing in Quantum Control literature \cite{Stefanatos03, Boscain05, Bonnard09, Lapert10}.

\section{FORMULATION OF THE PROBLEM IN TERMS OF OPTIMAL CONTROL}

The evolution of the wavefunction $\psi(t,x)$ of a particle in a
one-dimensional parabolic trapping potential with time-varying frequency
$\omega(t)$ is given by the Schr\"{o}dinger equation
\begin{equation}
\label{Schrodinger}i\hbar\frac{\partial\psi}{\partial t}=H(t)\psi=\left[  -\frac
{\hbar^{2}}{2m}\frac{\partial^{2}}{\partial x^{2}}+\frac{m\omega^{2}%
(t)}{2}x^{2}\right]  \psi,
\end{equation}
where $H(t)$ is the Hamiltonian operator of the system, $m$ is the particle mass and $\hbar$ is Planck's constant.
Consider the time evolution with initial frequency $\omega(0)=\omega_0$ at $t=0$ and final frequency $\omega(T)=\omega_T<\omega_0$ at the final time $T$. This corresponds to a temperature reduction (cooling) by a factor $\omega_T/\omega_0$. Frictionless cooling corresponds to a path $\omega(t)$ between these two values so that the populations of all the oscillator levels $n=0,1,2,\ldots$ at $t=T$ are equal to the ones at $t=0$. More explicitly, if
\begin{equation}
\psi(0,x)=\sum_{n=0}^{\infty}c_{n}(0)\Psi^{\omega_0}_{n}(x),\quad \psi(T,x)=\sum_{n=0}^{\infty}c_{n}(T)\Psi^{\omega_T}_{n}(x),
\end{equation}
where $\Psi^{\omega_0}_{n}(x),\Psi^{\omega_T}_{n}(x)$ are the eigenfunctions of the operators $H(0), H(T)$, respectively, then frictionless cooling corresponds to
\begin{equation}
\label{fric_cooling}
|c_n(0)|^2=|c_n(T)|^2,\quad n=0,1,2,\ldots
\end{equation}

It was shown in \cite{Chen10a}, using the Lewis-Riesenfeld invariant \cite{Lewis69}, that frictionless cooling is achieved when $\omega(t)$ satisfies the following Ermakov equation \cite{Ermakov}
\begin{equation}
\label{Ermakov}
\ddot{b}+\omega^2(t)b=\frac{\omega_0^2}{b^3},
\end{equation}
with boundary conditions
\begin{align}
\label{boundary1}b(0)&=1,\quad\dot{b}(0)=0,\quad\ddot{b}(0)=0,\\
\label{boundary2}b(T)&=\gamma,\quad\dot{b}(T)=0,\quad\ddot{b}(T)=0,
\end{align}
and $\gamma=\sqrt{\omega_0/\omega_T}>1$. Here $b(t)$ is a scaling dimensionless function describing the expansion of the wavefunction during the cooling process. When the above conditions are satisfied, the $n$th eigenstate of the initial oscillator at $t=0$ evolves following the ``expanding mode"
\begin{multline}
\label{singlemode}
\Psi_n(t,x)=\left(\frac{m\omega_0}{\pi\hbar}\right)^{1/4}\frac{\exp{\left[-i\left(n+\frac{1}{2}\right)\int_0^tdt'\frac{\omega_0}{b^2(t')}\right]}}{(2^nn!b)^{1/2}}\times
\\\exp{\left[i\frac{m}{2\hbar}\left(\frac{\dot{b}}{b}+\frac{i\omega_0}{b^2}\right)x^2\right]}
H_n\left[\left(\frac{m\omega_0}{\hbar}\right)^{1/2}\frac{x}{b}\right],
\end{multline}
where $H_n$ is the Hermite polynomial of degree $n$, and becomes eventually the $n$th eigenstate of the final trap at $t=T$, up to a global phase factor (independent of the spatial coordinate).
The instantaneous average energy $E_n(t)=\langle\Psi_n|H(t)|\Psi_n\rangle$ for the $n$th expanding mode is
\begin{equation}
\label{Inst_Energy}
E_n(t)=\frac{(2n+1)\hbar}{4\omega_0}\left[\dot{b}^2+\omega^2(t)b^2+\frac{\omega_0^2}{b^2}\right],
\end{equation}
and the time-averaged energy is
\begin{equation}
\label{Tim_Av_Energy}
\overline{E_n}=\frac{1}{T}\int_0^TE_n(t)dt.
\end{equation}
If we substitute (\ref{Inst_Energy}) in (\ref{Tim_Av_Energy}) and use the boundary conditions for $b$, we obtain
\begin{equation}
\label{Av_Energy}
\overline{E_n}=\frac{(2n+1)\hbar}{2\omega_0T}\int_0^T\left(\dot{b}^2+\frac{\omega_0^2}{b^2}\right)dt.
\end{equation}

We would like to find $\omega(t)$ satisfying the Ermakov equation (\ref{Ermakov}) and the corresponding boundary conditions (\ref{boundary1}) and (\ref{boundary2}), i.e., the sufficient conditions for frictionless cooling (\ref{fric_cooling}), while minimizing the time-averaged energy (\ref{Av_Energy}), for fixed final time $T$.
If we set
\begin{equation}
\label{definitions}
x_1=b,\quad x_2=\frac{\dot{b}}{\omega_0},\quad u(t)=\frac{\omega^2(t)}{\omega_0^2},
\end{equation}
and rescale time according to $t_{\mbox{new}}=\omega_0 t_{\mbox{old}}$, we obtain the following system of first order differential equations
\begin{align}
\label{system1}\dot{x}_1 &=x_2,\\
\label{system2}\dot{x}_2 &=-ux_1+\frac{1}{x_1^3},
\end{align}
equivalent to the Ermakov equation (\ref{Ermakov}), and
\begin{equation}
\label{Short_Av_Energy}
\overline{E_n}=\frac{(2n+1)\hbar\omega_0}{T}J,
\end{equation}
where
\begin{equation}
\label{cost}
J=\frac{1}{2}\int_0^T\left(x_2^2+\frac{1}{x_1^2}\right)dt.
\end{equation}

The corresponding optimal control problem is: Given the system (\ref{system1}), (\ref{system2}), with initial condition $(x_1(0),x_2(0))=(1,0)$ and final condition $(x_1(T),x_2(T))=(\gamma,0)$, find the control $u(t), 0\leq t \leq T$, $T$ fixed, with $u(0)=1, u(T)=1/\gamma^4$, that minimizes the cost $J$ given in (\ref{cost}). The boundary conditions on the state variables $(x_1,x_2)$ are equivalent to those for $b$ and $\dot{b}$, while the boundary conditions on the control variable $u$ lead to the corresponding conditions for $\ddot{b}$. Note that the possibility $\omega^2(t)<0$ (expulsive parabolic potential) for some time intervals is permitted \cite{Chen10a}. Also, for fixed final time $T$, the minimization of cost $J$ corresponds to minimizing the time-averaged energy (\ref{Short_Av_Energy}).

If we omit the boundary conditions on $u(t)$ and solve the corresponding
problem, the minimum cost that we find is a lower bound of the minimum cost
for the full problem, where these conditions are on. In the following we solve the relaxed problem, while the study of the full case will be the subject of a subsequent publication.

\newtheorem{problem}{Problem}
\begin{problem}
\label{problem1}
Given the system (\ref{system1}), (\ref{system2}) with initial condition $(x_1(0),x_2(0))=(1,0)$ and final condition $(x_1(T),x_2(T))=(\gamma,0), \gamma>1$, find the control $u(t), 0\leq t \leq T$, $T$ fixed, that minimizes the cost $J$ given in (\ref{cost}).
\end{problem}


\section{SINGULAR SOLUTION}

The control Hamiltonian is given by
\begin{equation}
\label{Hamiltonian}
H_c=\lambda_0\frac{1}{2}\left(x_2^2+\frac{1}{x_1^2}\right)+\lambda_1x_2+\lambda_2\left(-ux_1+\frac{1}{x_1^3}\right),
\end{equation}
where, according to Maximum Principle \cite{Pontryagin}, $\lambda_0\leq 0$ is a constant and the adjoint variables satisfy the equations
\begin{align}
\label{lambda_1}\dot{\lambda}_1 &=-\frac{\partial H_c}{\partial x_1}=\lambda_0\frac{1}{x_1^3}+\lambda_2\left(u+\frac{3}{x_1^4}\right),\\
\label{lambda_2}\dot{\lambda}_2 &=-\frac{\partial H_c}{\partial x_2}=-\lambda_0x_2-\lambda_1.
\end{align}

Observe that $H_c$ is linear in $u$. The switching function is $\Phi=-\lambda_2$ (note that $x_1>0$ due to the repulsive force $1/x_1^3$ at $x_1=0$). Singular arcs are encountered when $\lambda_2=0$ for some finite time interval. Note that if also $\lambda_0=0$ then from (\ref{lambda_2}) it is $\lambda_1=0$, which is forbidden since Maximum Principle requires $(\lambda_0,\lambda_1,\lambda_2)\neq 0$ \cite{Pontryagin}. So $\lambda_0\neq0$ on a singular arc and we set $\lambda_0=-1$. For $\lambda_2=0$, (\ref{lambda_1}), (\ref{lambda_2}) reduce to $\dot{\lambda}_1=-1/x_1^3, \lambda_1=x_2$, respectively. Setting $\lambda_1=x_2$  in (\ref{Hamiltonian}) we obtain
\begin{equation}
\label{Sing_Ham}
H_c=\frac{1}{2}\left(x_2^2-\frac{1}{x_1^2}\right).
\end{equation}
Since the system (\ref{system1}), (\ref{system2}) is autonomous, the control Hamiltonian is constant. The equation
\begin{equation}
\label{Sing_Arcs}
x_2^2-\frac{1}{x_1^2}=c
\end{equation}
represents a one-parameter family of singular arcs in state space. Combining the adjoint equations we find $\dot{x}_2=-1/x_1^3$, and using this in (\ref{system2}) we obtain the expression for the control on the singular arc
\begin{equation}
\label{Sing_Control}
u_s=\frac{2}{x_1^4}
\end{equation}

\subsection{Unbounded Control}

\begin{figure}[t]
\centering
\includegraphics[width=0.6\linewidth]{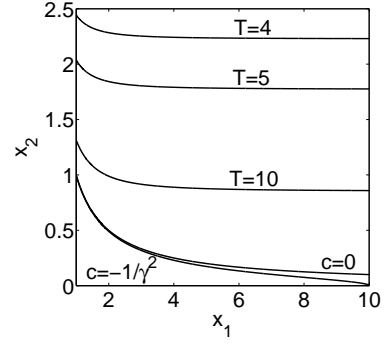}
\caption{Singular trajectories for various values of the final time $T$ when $\gamma=10$ and the control is unbounded. The trajectories characterized by the constant values of the control Hamiltonian $c=0$, $c=-1/\gamma^2$, correspond to $T=(\gamma^2-1)/2$, $T=\gamma\sqrt{\gamma^2-1}$, respectively. The initial impulse, driving the starting state $(1,0)$ to the singular arc is negative, while the final impulse, driving the state from the singular arc to the target point $(\gamma,0)$, is positive.}
\label{fig:unbound_traj1}
\end{figure}

\begin{figure}[t]
\centering
\includegraphics[width=0.6\linewidth]{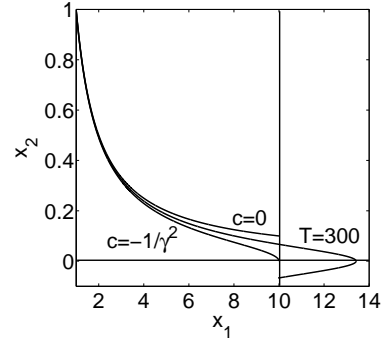}
\caption{Singular trajectory for $T>\gamma\sqrt{\gamma^2-1}$. Note that final impulse is negative.}
\label{fig:unbound_traj2}
\end{figure}

When the control is unbounded, the state can be shifted along the lines of constant $x_1$ instantaneously by the use of Dirac impulses in $u$. Such controls have no effect on the performance index $J$ directly, since $u(t)$ does not enter the cost function. A typical extremal solution involves an initial impulse to move the state to the singular arc (at $t=0^{+}$), then motion along the singular arc until the line $x_1=\gamma$ is reached, then another impulse to move the state to the target point $(\gamma,0)$. The fact that the state must arrive at the target point at $t=T$ determines the constant $c$ in (\ref{Sing_Arcs}), which picks out the particular singular arc in the one-parameter family of possible arcs. Using (\ref{system1}) and (\ref{Sing_Arcs}), one can find the time evolution of $x_1$ along the singular arc
\begin{equation}
\label{x1}
x_1(t)=\sqrt{(B^2-T^2)\left(\frac{t}{T}\right)^2+2B\left(\frac{t}{T}\right)+1},
\end{equation}
where $B=\sqrt{\gamma^2+T^2}-1$, in agreement with the result obtained in \cite{Chen10b} using the Euler-Lagrange equation.
The constant defining the singular arc in (\ref{Sing_Arcs}) is $c=(B/T)^2-1$. The initial impulse moving the state to the singular arc is negative. Now
\begin{equation}
x_2(T^-)=\dot{x}_1(T)=\frac{\gamma^2-\sqrt{\gamma^2+T^2}}{\gamma T},
\end{equation}
so $x_2(T^-)>0$ for $T<\gamma\sqrt{\gamma^2-1}$ and $x_2(T^-)<0$ for $T>\gamma\sqrt{\gamma^2-1}$. The final impulse that drives the state to $x_2(T^+)=0$ is positive in the first case and negative in the second. Characteristic trajectories in state space for these two cases are plotted in Figs. \ref{fig:unbound_traj1} and \ref{fig:unbound_traj2}, respectively. It can be easily shown that $dc/dT<0$ for $T<\gamma\sqrt{\gamma^2-1}$ and $dc/dT>0$ for $T>\gamma\sqrt{\gamma^2-1}$. The corresponding minimum value is $c_{\mbox{min}}=-1/\gamma^2$ and $c\in[-1/\gamma^2,\infty)$, with $c\rightarrow \infty$ for $T\rightarrow 0$ and $c\rightarrow 0$ for $T\rightarrow \infty$. Note that as the final time $T$ varies, the minimum cost $J$ is obtained when $H_c=0\Rightarrow c=0\Rightarrow T=(\gamma^2-1)/2$. This case is equivalent to finding the minimum $J$ with the final time $T$ free. Returning to the fixed time case, we are interested in the time-averaged energy $\overline{E}=2J/T$, in units of $(n+1/2)\hbar\omega_0$, which is given by the following expression
\begin{equation}
\label{Extr_Energy1}
\overline{E}=\frac{1}{T^2}\left[\gamma^2+1-2\sqrt{\gamma^2+T^2}+2T\ln{\left(\frac{T+\sqrt{\gamma^2+T^2}}{\gamma}\right)}\right]
\end{equation}
and is plotted versus $T$ in Fig. \ref{fig:unbound_energy}. It is shown in Appendix \ref{appendix1} that $d\overline{E}/dT<0$, so the time-averaged energy is a decreasing function of the final time $T$.

\begin{figure}[t]
\centering
\includegraphics[width=0.6\linewidth]{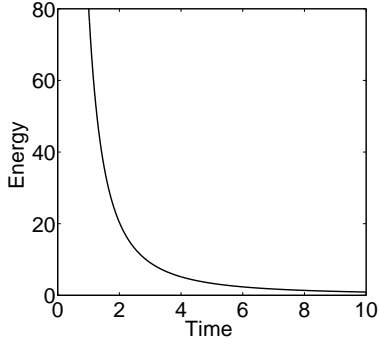}
\caption{Time-averaged energy $\overline{E}$, in units of $(n+1/2)\hbar\omega_0$, plotted as a function of the final time $T$, for $\gamma=10$ and unbounded control. Note that $d\overline{E}/dT<0$.}
\label{fig:unbound_energy}
\end{figure}

There is an alternative way to obtain the above results and the optimality of the singular solution. The idea is to define a ``dual" time-optimal problem and use the results regarding singular trajectories of such problems \cite{Bonnard03}. With this in mind, we augment the system (\ref{system1}), (\ref{system2}) with a state variable $x_3$ satisfying the differential equation
\begin{equation}
\label{system3}
\dot{x}_3=x_2^2+\frac{1}{x_1^2}-\overline{E}
\end{equation}
and the boundary conditions $x_3(0)=x_3(T)=0$, where $\overline{E}>0$ is fixed. We now seek a control that drives the state of the augmented system from $x(0)=(1,0,0)$ to $x(0)=(\gamma,0,0), \gamma>1$, in minimum time $T$. Observe that the extra state is defined such that the equality
\begin{equation}
\label{fixed_energy}
\frac{1}{T}\int_0^T\left(x_2^2+\frac{1}{x_1^2}\right)dt=\overline{E}
\end{equation}
holds, in other words the time-averaged energy is now fixed while the final time is free. The corresponding affine control system is $\dot{x}=f(x)+ug(x)$, where the vector fields $f,g$ are
\begin{equation}
\label{vector_fields}
f=(x_2, \frac{1}{x_1^3}, x_2^2+\frac{1}{x_1^2}-\overline{E})^T,\quad g=(0, -x_1, 0)^T.\nonumber
\end{equation}
The next step is to calculate the following determinants
\begin{align}
D&=\mbox{det}(g,[f,g],[g,[f,g]])=-2x_1^4,\nonumber\\
D'&=\mbox{det}(g,[f,g],[f,[f,g]])=4,\nonumber\\
D''&=\mbox{det}(g,[f,g],f)=1-x_1^2(\overline{E}+x_2^2),\nonumber
\end{align}
where the Lie-bracket is $[f,g]=(\partial g/\partial x)f-(\partial f/\partial x)g$. The singular control for the augmented time-optimal problem is given by $u_s=-D'/D=2/x_1^4$ \cite{Bonnard03}, and is the same as in (\ref{Sing_Control}). As a consequence, the projection of the singular arc on the $x_1x_2$-plane is given by (\ref{Sing_Arcs}) and $T$ is determined implicitly by (\ref{Extr_Energy1}), since $\overline{E}$ is fixed. We now show that this solution corresponds to a minimum. We first show that the singular arc is hyperbolic. The criterion for this is \cite{Bonnard03}
\begin{equation}
DD''>0\Leftrightarrow \overline{E}+x_2^2-\frac{1}{x_1^2}=\overline{E}+c>0.
\end{equation}
Using (\ref{Extr_Energy1}) and the expression for $c$ given above we find $\overline{E}+c=2f(T)/T^2$, where the function $f(T)>0$ is given in Appendix \ref{appendix1}. So the singular arc is hyperbolic. In this case, we know that it corresponds to the minimum time solution, up to the first conjugate time $t_{1c}$ \cite{Bonnard03}. In order to find $t_{1c}$ we introduce the vector field $S=f+u_sg$, that generates the singular trajectory, and the Jacobi field $V$, which is the solution of the variational equation
\begin{equation}
\dot{\delta x}=\frac{\partial S}{\partial x}\delta x,\nonumber
\end{equation}
with the initial condition $\delta x(0)=g(0)=(0,-1,0)^T$. The first conjugate time is the first time such that the fields $V,g$ are collinear, i.e., $V(t_{1c})\parallel (0,1,0)^T$. The first two equations of the variational system are
\begin{align}
\dot{\delta x_1}&=\delta x_2,\nonumber\\
\dot{\delta x_2}&=\frac{3}{x_1^4}\delta x_1,\nonumber
\end{align}
with the initial conditions $\delta x_1(0)=0, \delta x_2(0)=-1$. It is not hard to see that $\delta x_1(t)<0$ for $t>0$, so $V$ cannot be aligned with $g$ and there is no conjugate point. As a result, the singular solution is time-optimal for the augmented system with fixed time-averaged energy. This elucidates the conjugate relation between energy and time.

\subsection{Bounded Control}

In practice there are limits on the control amplitude that set a tighter lower bound on the transfer time. To fix ideas, we consider the case $|u|\leq 1$ for the two-dimensional system (\ref{system1}), (\ref{system2}). We have shown in our previous work \cite{Stefanatos11} that in this case the minimum necessary time to transfer the initial state $(1,0)$ to the target state $(\gamma,0), \gamma>1$, is obtained following the bang-bang strategy
\begin{equation}
\label{bang-bang}
u(t)=\left\{\begin{array}{cl} -1, & t<T_1 \\1, & T_1<t\leq T_1+T_2\end{array}\right..
\end{equation}
From (\ref{system1}), (\ref{system2}) we find the trajectory for $t<T_1$
\begin{equation}
\label{trajectory1}
x_2^2-x_1^2+\frac{1}{x_1^2}=0,
\end{equation}
and the trajectory for $T_1<t\leq T_1+T_2$
\begin{equation}
\label{trajectory2}
x_2^2+x_1^2+\frac{1}{x_1^2}=\gamma^2+\frac{1}{\gamma^2}.
\end{equation}
The $x_1$ coordinate of their common point at $t=T_1$ is
\begin{equation}
\label{joint}
x_1(T_1)=\sqrt{\frac{\gamma^4+1}{2\gamma^2}}
\end{equation}
Integrating (\ref{system1}), (\ref{system2})  from the initial point up to the common point, for duration $t=T_1$ and $u=-1$, we find
\begin{equation}
\label{T1}
x_1(T_1)=\sqrt{\cosh(2T_1)}.
\end{equation}
Integrating (\ref{system1}), (\ref{system2}) from the final point back to the common point, for duration $t=T_2$ and $u=1$, we find
\begin{equation}
\label{T2}
x_1(T_1)=\sqrt{\frac{1}{2}\left(\gamma^2+\frac{1}{\gamma^2}\right)+\frac{1}{2}\left(\gamma^2-\frac{1}{\gamma^2}\right)\cos(2T_2)}
\end{equation}
From (\ref{joint}), (\ref{T1}) and (\ref{T2}) we obtain
\begin{equation}
\label{MinTime}
T_{\mbox{min}}=T_1+T_2=\frac{1}{2}\cosh^{-1}{\left(\frac{\gamma^4+1}{2\gamma^2}\right)}+\frac{\pi}{4}
\end{equation}

\begin{figure}[t]
\centering
\includegraphics[width=0.6\linewidth]{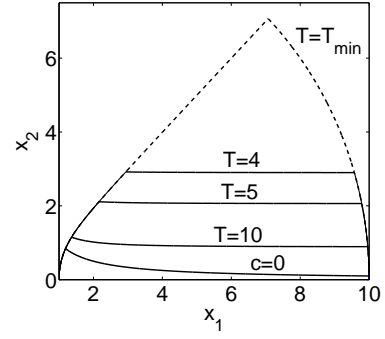}
\caption{Trajectories corresponding to the (bounded) control policy (\ref{Sing_Con1}), for various values of the final time $T\geq T_{\mbox{min}}$ and $\gamma=10$.}
\label{fig:bound_traj}
\end{figure}

We see that in the presence of control bounds, problem \ref{problem1} is meaningful for final times  $T\geq T_{\mbox{min}}$.
For simplicity, we will restrict the rest of the discussion to the bang-singular-bang case
\begin{equation}
\label{Sing_Con1}
u(t)=\left\{\begin{array}{cl} -1, & t<T'_1 \\ 2/x_1^4, &T'_1<t\leq T'_1+T'_2 \\1, & T'_1+T'_2<t\leq T'_1+T'_2+T'_3\end{array}\right..
\end{equation}
Under this control policy, the initial and final arcs (\ref{trajectory1}) and (\ref{trajectory2}) are joined by a singular arc (\ref{Sing_Arcs}), as it is depicted in Fig. \ref{fig:bound_traj}. We expect that at least for the most interesting case, where $T$ is close to $T_{\mbox{min}}$, the pulse sequence (\ref{Sing_Con1}) is the optimal one. Note that the joint point for $T=T_{\mbox{min}}$ lies above the singular arc $x_2=1/x_1$ corresponding to $c=0$, as shown in Fig. \ref{fig:bound_traj}, when $\gamma^8-6\gamma^4+1>0\Rightarrow\gamma>\sqrt[4]{3+2\sqrt{2}}\approx 1.554$. Throughout the text we use the value $\gamma=10$.

\begin{figure}[t]
\centering
\includegraphics[width=0.6\linewidth]{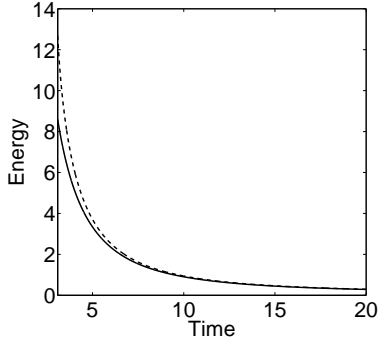}
\caption{Time-averaged energy corresponding to the (bounded) control policy (\ref{Sing_Con1}), as a function of the final time $T\geq T_{\mbox{min}}$, for $\gamma=10$ (dashed line). The corresponding energy for the case of unbounded control and the same target point is also plotted (solid line).}
\label{fig:bound_energy}
\end{figure}

Using (\ref{trajectory1}), (\ref{Sing_Arcs}) we find for the first junction point
\begin{equation}
\label{x1T1}
x_1(T'_1)=\sqrt{\frac{c+\sqrt{c^2+8}}{2}},
\end{equation}
while from (\ref{trajectory2}), (\ref{Sing_Arcs}) we find for the second junction point
\begin{equation}
\label{x1T2}
x_1(T'_1+T'_2)=\sqrt{\frac{-c\gamma^2+\gamma^4+1+\sqrt{(c\gamma^2-\gamma^4-1)^2-8\gamma^4}}{2\gamma^2}},
\end{equation}
If we limit our analysis to the cases where the constant determining the singular arc satisfies $c>0$, which also correspond to shorter transfer times, then
\begin{equation}
u_s(t)=\frac{2}{x_1^4(t)}\leq\frac{2}{x_1^4(T'_1)}=\frac{8}{(c+\sqrt{c^2+8})^2}< 1,
\end{equation}
so the control remains within the allowed bounds along the singular arc.
By integrating the equations of motion we obtain the time spent on each arc
\begin{align}
\label{time1}T'_1(c) &=\frac{1}{2}\cosh^{-1}{x_1^2(T'_1)},\\
\label{time2}T'_2(c) &=\frac{x_1^2(T'_1+T'_2)-x_1^2(T'_1)}{\sqrt{1+cx_1^2(T'_1+T'_2)}+\sqrt{1+cx_1^2(T'_1)}},\\
\label{time3}T'_3(c) &=\frac{1}{2}\cos^{-1}\left(\frac{2\gamma^2x_1^2(T'_1+T'_2)-\gamma^4-1}{\gamma^4-1}\right).
\end{align}
The constant $c$ which determines the singular arc can be found from the transcendental equation
\begin{equation}
\label{transcendental1}
T'_1(c)+T'_2(c)+T'_3(c)=T.
\end{equation}
Having determined $c$, we can calculate the time-averaged energy corresponding to the control sequence (\ref{Sing_Con1}) for a specific final time. The result is plotted in Fig. \ref{fig:bound_energy}, along with the time-averaged energy that derived previously for the unbounded control case. As expected, when the control is bounded, a larger amount of energy is necessary to achieve the same transfer for times close to $T_{\mbox{min}}$. Equivalently, for a fixed level of time-averaged energy, a longer time is needed to reach the same target point.

\section{CONCLUSION AND FUTURE WORK}

In this paper we examined from a control-theoretic viewpoint the frictionless cooling of atoms trapped in a harmonic potential, while minimizing the transient energy of the system, a problem that was first considered in \cite{Chen10b}. We showed that the transient energy is minimized by a singular control. We also showed that this singular control is the time-minimal solution for a ``dual" problem, where the energy is held fixed, highlighting the conjugate relation between transient energy and minimum time. In addition, we examined briefly how the solution is modified when there are bounds on the control.

Possible future work could include the detailed examination of the bounded control case, that was only slightly touched here, as well as the incorporation in the analysis of additional restrictions on the controls reflecting experimental limitations \cite{Schaff10}. The complexity of the resulting optimal control problems, which may increase the difficulty of the analytical study, can be overcome by using a powerful state of the art numerical optimization method based on pseudospectral approximations \cite{Li09, Li_PNAS11}.

The results presented here can be immediately extended to the frictionless cooling of a two-dimensional Bose-Einstein condensate confined in a parabolic trapping potential \cite{Muga09, Schaff_EPL11}. The above techniques are not restricted to atom cooling but are applicable to areas as diverse as adiabatic quantum computing \cite{Aharonov07} and finite time thermodynamic processes \cite{Salamon09}.

\section{ACKNOWLEDGMENTS}

The authors would like to thank Profs. J. G. Muga and H. Schaettler for valuable comments.

\addtolength{\textheight}{-15.5cm}

\appendices
\section{}

\label{appendix1}

We show that $d\overline{E}/dT<0$. From (\ref{Extr_Energy1}) we obtain $d\overline{E}/dT=-2f(T)/T^3$, where
\begin{equation}
f(T)=\gamma^2+1-2\sqrt{\gamma^2+T^2}+T\ln{\frac{T+\sqrt{\gamma^2+T^2}}{\gamma}}.
\end{equation}
It is $f(0)=(\gamma-1)^2>0$ and $df/dT=g(T)$, where
\begin{equation}
g(T)=\ln{\frac{T+\sqrt{\gamma^2+T^2}}{\gamma}}-\frac{T}{\sqrt{\gamma^2+T^2}}.
\end{equation}
It is $g(0)=0$ and
\begin{equation}
\frac{dg}{dT}=T^2(\gamma^2+T^2)^{-3/2}\geq 0.
\end{equation}
So $g(T)\geq 0\Rightarrow f(T)>0\Rightarrow dE_n/dT<0$.

\end{document}